\documentclass[12pt,twoside]{article}

\usepackage[english]{babel}
\usepackage{umj_1}           

\usepackage{graphicx, srcltx}   
\usepackage{cite}               

\usepackage{amsmath}            
\usepackage{amsfonts,amssymb}   
\usepackage{srcltx}             

\newcommand{\Leftclt}{Roman A. Veprintsev}
\newcommand{\Rightclt}{On the asymptotic behavior of generalized Gegenbauer polynomials}
\usepackage{authblk}

\begin{document}

\begin{Titul}
{\large \bf  ON THE ASYMPTOTIC BEHAVIOR\\ OF THE MAXIMUM ABSOLUTE VALUE\\[0.2em] OF GENERALIZED GEGENBAUER POLYNOMIALS}\\[3ex]
{{\bf Roman~A.~Veprintsev} \\[5ex]}
\end{Titul}

\begin{Anot}
{\bf Abstract.} Using well-known facts on Jacobi polynomials, we derive some asymptotic estimates for the maximum absolute value of generalized Gegenbauer polynomials.


\end{Anot}


\section*{Introduction and main result}

In this section, we introduce some classes of orthogonal polynomials on $[-1,1]$, including the so-called generalized Gegenbauer polynomials. For a background and more details on the orthogonal polynomials, the reader is referred to \cite{dai_xu_book_approximation_theory_2013,dunkl_xu_book_orthogonal_polynomials_2014,andrews_askey_roy_book_special_functions_1999,szego_book_orthogonal_polynomials_1975}. Here, we also formulate the main result of the publication.

Let $\alpha,\,\beta>-1$. The Jacobi polynomials, denoted by $P_n^{(\alpha,\beta)}(\cdot)$, where $n=0,1,\ldots$, are orthogonal with respect to the Jacobi weight function $w_{\alpha,\beta}(t)=(1-t)^\alpha(1+t)^\beta$ on $[-1,1]$, namely,
\begin{equation*}
\int\nolimits_{-1}^1 P_n^{(\alpha,\beta)}(t)\, P_m^{(\alpha,\beta)}(t)\, w_{\alpha,\beta}(t)\,dt=\begin{cases}\dfrac{2^{\alpha+\beta+1}\Gamma(n+\alpha+1)\Gamma(n+\beta+1)}{(2n+\alpha+\beta+1)\Gamma(n+1)\Gamma(n+\alpha+\beta+1)},&n=m,\\
0,&n\not=m.
\end{cases}
\end{equation*}
Here, as usual, $\Gamma$ is the gamma function.

For $\lambda>-\frac{1}{2}$, $\mu\geq0$, and $n=0,1,\ldots$, the generalized Gegenbauer polynomials $C_n^{(\lambda,\mu)}(\cdot)$ are defined by
\begin{equation*}\label{coefficients_for_generalized_Gegenbauer_polynomials}
\begin{array}{ll}
C_{2n}^{(\lambda,\mu)}(t)=a_{2n}^{(\lambda,\mu)}P_n^{(\lambda-1/2,\mu-1/2)}(2t^2-1), & a_{2n}^{(\lambda,\mu)}=\dfrac{(\lambda+\mu)_n}{(\mu+\frac{1}{2})_n},\\[1.0em]
C_{2n+1}^{(\lambda,\mu)}(t)=a_{2n+1}^{(\lambda,\mu)}\,t P_n^{(\lambda-1/2,\mu+1/2)}(2t^2-1),\quad & a_{2n+1}^{(\lambda,\mu)}=\dfrac{(\lambda+\mu)_{n+1}}{(\mu+\frac{1}{2})_{n+1}},
\end{array}
\end{equation*}
where $(\lambda)_n$ denotes the Pochhammer symbol given by
\begin{equation*}
(\lambda)_0=1,\quad (\lambda)_n=\lambda(\lambda+1)\cdots(\lambda+n-1)\quad\text{ for}\quad n=1,2,\ldots.
\end{equation*}
They are orthogonal with respect to the weight function
\[v_{\lambda,\mu}(t)=|t|^{2\mu}(1-t^2)^{\lambda-1/2},\quad t\in[-1,1].\]
For $\mu=0$, these polynomials, denoted by $C_n^{\lambda}(\cdot)$, are called the Gegenbauer polynomials:
\begin{equation*}
C_n^{\lambda}(t)=C_n^{(\lambda,0)}(t)=\frac{(2\lambda)_n}{(\lambda+\frac{1}{2})_n} P_n^{(\lambda-1/2,\lambda-1/2)}(t).
\end{equation*}

For $\lambda>-\frac{1}{2}$, $\mu>0$, and $n=0,1,\ldots$, we have the following connection:
\begin{equation*}
C_n^{(\lambda,\mu)}(t)=c_\mu\int\nolimits_{-1}^1 C_n^{\lambda+\mu}(tx)(1+x)(1-x^2)^{\mu-1}\,dx,\quad c_\mu^{-1}=2\int\nolimits_0^1 (1-x^2)^{\mu-1}\,dx.
\end{equation*}

Denote by $\bigl\{\widetilde{C}_n^{(\lambda,\mu)}(\cdot)\bigr\}_{n=0}^{\infty}$ the sequence of orthonormal generalized Gegenbauer polynomials. It is easily verified that these polynomials are given by the following formulae:
\begin{equation}\label{coefficients_for_orthonormal_generalized_Gegenbauer_polynomials}
\begin{split}
&\widetilde{C}_{2n}^{(\lambda,\mu)}(t)=\widetilde{a}_{2n}^{\,(\lambda,\mu)}P_n^{(\lambda-1/2,\mu-1/2)}(2t^2-1),\\ &\widetilde{a}_{2n}^{\,(\lambda,\mu)}=\Biggl(\dfrac{(2n+\lambda+\mu)\Gamma(n+1)\Gamma(n+\lambda+\mu)}{\Gamma(n+\lambda+\frac{1}{2})\Gamma(n+\mu+\frac{1}{2})}\Biggr)^{1/2},\\[0.4em]
&\widetilde{C}_{2n+1}^{(\lambda,\mu)}(t)=\widetilde{a}_{2n+1}^{\,(\lambda,\mu)}\,t P_n^{(\lambda-1/2,\mu+1/2)}(2t^2-1),\\ &\widetilde{a}_{2n+1}^{\,(\lambda,\mu)}=\Biggl(\dfrac{(2n+\lambda+\mu+1)\Gamma(n+1)\Gamma(n+\lambda+\mu+1)}{\Gamma(n+\lambda+\frac{1}{2})\Gamma(n+\mu+\frac{3}{2})}\Biggr)^{1/2}.
\end{split}
\end{equation}

Throughout the paper we use the following asymptotic notation: $f(n)\lesssim g(n)$, $n\to\infty,$ or equivalently $g(n)\gtrsim f(n)$, $n\to\infty$, means that there exist a positive constant $C$ and a positive integer $n_0$ such that $0\leq f(n)\leq C g(n)$ for all $n\geq n_0$ (asymptotic upper bound); if there exist positive constants $C_1$, $C_2$, and a positive integer $n_0$ such that $0\leq C_1 g(n)\leq f(n)\leq C_2 g(n)$ for all $n\geq n_0$, then we write $f(n)\asymp g(n)$, $n\to\infty$ (asymptotic tight bound).

To simplify the writing, we will omit ``$n\to\infty$'' in the asymptotic notation.

It follows directly from Stirling's asymptotic formula that
\begin{equation*}
\frac{\Gamma(n+\alpha)}{\Gamma(n+\beta)}\asymp n^{\alpha-\beta}
\end{equation*}
for arbitrary real numbers $\alpha$ and $\beta$. For appropriate values of $q$, $n$ the formula $\frac{\Gamma(q+n)}{\Gamma(q)}=(q)_n$ holds. Thus, for $q>0$,
\begin{equation*}
\frac{(q)_n}{n!}=\frac{(q)_n}{\Gamma(n+1)}\asymp n^{q-1}.
\end{equation*}
Hence, by \eqref{coefficients_for_orthonormal_generalized_Gegenbauer_polynomials},
\begin{equation}\label{asymptotic_for_coefficients_for_orthonormal_generalized_Gegenbauer_polynomials}
\widetilde{a}_{2n}^{\,(\lambda,\mu)}\asymp n^{1/2},\quad \widetilde{a}_{2n+1}^{\,(\lambda,\mu)}\asymp n^{1/2}.
\end{equation}

Define the uniform norm of a continuous function $f$ on $[-1,1]$ by
\begin{equation*}
\|f\|_{\infty}=\max\limits_{-1\leq t\leq 1} |f(t)|.
\end{equation*}
The maximum of two real numbers $x$ and $y$ is denoted by $\max(x,y)$.

Now we can formulate the main result.

\begin{teoen}\label{main_result_theorem}
Let $\lambda>-\frac{1}{2}$, $\mu>0$. Then 
\begin{equation*}
\bigl\|\widetilde{C}_{n}^{(\lambda,\mu)}\bigr\|_\infty\asymp n^{\max(\lambda,\mu)}.
\end{equation*}
\end{teoen}

The proof of the above theorem is contained in Section \ref{section_proof_of_the_main_result}. In the next section, we give some important properties of Jacobi polynomials needed for our purpose.

\section{Some facts on Jacobi polynomials}

For $\alpha,\,\beta>-1$, the following formulae are valid \cite[\S\,2.3, Corollary~3.2]{osilenker_book_orthogonal_polynomials_1999}
\begin{equation}\label{connection_with_plus_minus}
P_n^{(\alpha,\beta)}(-t)=(-1)^n P_n^{(\beta,\alpha)}(t),
\end{equation}
\begin{equation}\label{values_at_end_points}
P_n^{(\alpha,\beta)}(1)=\frac{(\alpha+1)_n}{n!}\asymp n^\alpha,\quad |P_n^{(\alpha,\beta)}(-1)|=\frac{(\beta+1)_n}{n!}\asymp n^\beta.
\end{equation}

We have \cite[Theorem~7.32.1]{szego_book_orthogonal_polynomials_1975}
\begin{equation}\label{asymptotic_for_maximum_of_Jacobi_polynomials}
\bigl\|P_n^{(\alpha,\beta)}\bigr\|_\infty=\left\{\begin{array}{lrclll}
P_n^{(\alpha,\beta)}(1),&\alpha\!\!\!&\geq\beta,&&\alpha\geq&\!\!\!-\frac{1}{2},\\[0.4em]
|P_n^{(\alpha,\beta)}(-1)|,&\alpha\!\!\!&\leq\beta,&&\beta\geq&\!\!\!-\frac{1}{2}
\end{array}\right.
\end{equation}
and
\begin{equation}\label{asymptotic_of_Jacobi_polynomials_on_half-segment}
\max\limits_{0\leq t\leq 1} |P_n^{(\alpha,\beta)}(t)|\lesssim n^{\max(\alpha,-1/2)}.
\end{equation}

It follows from Theorem~8.1.1 in \cite{szego_book_orthogonal_polynomials_1975} that, for $\alpha>-\frac{1}{2}$,
\begin{equation}\label{asymptotic_of_Jacobi_polynomials_on_the_sequence_of_special_points}
|P_n^{(\alpha,\beta)}(\cos n^{-1})|\asymp n^\alpha.
\end{equation}

It is known (see, for example, \cite[Theorem~7.32.2]{szego_book_orthogonal_polynomials_1975}) that
\begin{equation}\label{asymptotic_of_Jacobi_polynomials_on_dynamic_parts_of_segment}
\begin{split}
&\max\limits_{\frac{1}{n}\leq\theta\leq\frac{\pi}{2}} \Bigl\{\theta^{\,\alpha+\frac{1}{2}} |P_n^{(\alpha,\beta)}(\cos\theta)|\Bigr\}\lesssim n^{-\frac{1}{2}},\\
&\max\limits_{0\leq\theta\leq\frac{1}{n}} |P_n^{(\alpha,\beta)}(\cos\theta)|\lesssim n^{\alpha}.
\end{split}
\end{equation}

\begin{lemen}\label{main_result_lemma}
Let $\alpha>\frac{1}{2}$. Then
\begin{equation*}
\begin{split}
&\max\limits_{0\leq\theta\leq\frac{\pi}{2}} \Bigl\{\sin\frac{\theta}{2}\, |P_n^{(\alpha,\beta)}(\cos\theta)|\Bigr\}\lesssim n^{\alpha-1},
\\
&\max\limits_{\frac{\pi}{2}\leq\theta\leq \pi} \Bigl\{\sin\frac{\theta}{2}\, |P_n^{(\alpha,\beta)}(\cos\theta)|\Bigr\}\lesssim n^{\max(\beta,-1/2)}.
\end{split}
\end{equation*}

\end{lemen}

\proofen Using the assumption $\alpha>\frac{1}{2}$ and the asymptotic equality $\sin\frac{n^{-1}}{2}\asymp n^{-1}$, it can be easily seen that, for sufficiently large $n$,
\begin{equation*}
\max\limits_{\frac{1}{n}\leq\theta\leq\frac{\pi}{2}} \Bigl\{\theta^{\,-\alpha-\frac{1}{2}} \sin\frac{\theta}{2}\Bigr\}=n^{\alpha+\frac{1}{2}} \sin\frac{n^{-1}}{2}\asymp n^{\alpha-1/2}.
\end{equation*}
Hence, by \eqref{asymptotic_of_Jacobi_polynomials_on_dynamic_parts_of_segment}, we obtain
\begin{equation*}
\max\limits_{0\leq\theta\leq\frac{1}{n}} \Bigl\{\sin\frac{\theta}{2}\,|P_n^{(\alpha,\beta)}(\cos\theta)|\Bigr\}\lesssim n^{\alpha-1},
\end{equation*}
\begin{equation*}
\begin{split}
\max\limits_{\frac{1}{n}\leq\theta\leq\frac{\pi}{2}} &\Bigl\{\sin\frac{\theta}{2}\, |P_n^{(\alpha,\beta)}(\cos\theta)|\Bigr\}\leq\\
&\leq\max\limits_{\frac{1}{n}\leq\theta\leq\frac{\pi}{2}} \Bigl\{\theta^{\,-\alpha-\frac{1}{2}} \sin\frac{\theta}{2}\Bigr\}\, \cdot \, \max\limits_{\frac{1}{n}\leq\theta\leq\frac{\pi}{2}} \Bigl\{\theta^{\,\alpha+\frac{1}{2}} |P_n^{(\alpha,\beta)}(\cos\theta)|\Bigr\}\lesssim n^{\alpha-1}.
\end{split}
\end{equation*}
So, we get the desired estimate on $\bigl[0,\frac{\pi}{2}\bigr]$.

The asymptotic upper bound on $\bigr[\frac{\pi}{2},\pi\bigl]$ follows from \eqref{connection_with_plus_minus}, \eqref{asymptotic_of_Jacobi_polynomials_on_half-segment}. Indeed,
\begin{equation*}
\begin{split}
\max\limits_{\frac{\pi}{2}\leq\theta\leq \pi} &\Bigl\{\sin\frac{\theta}{2}\, |P_n^{(\alpha,\beta)}(\cos\theta)|\Bigr\}\asymp \max\limits_{\frac{\pi}{2}\leq\theta\leq \pi} |P_n^{(\alpha,\beta)}(\cos\theta)|=\max\limits_{0\leq\theta\leq \frac{\pi}{2}} |P_n^{(\alpha,\beta)}(\cos(\pi-\theta))|=\\
&=\max\limits_{0\leq\theta\leq \frac{\pi}{2}} |P_n^{(\alpha,\beta)}(-\cos\theta)|=\max\limits_{0\leq\theta\leq \frac{\pi}{2}} |P_n^{(\beta,\alpha)}(\cos\theta)|\lesssim n^{\max(\beta,-1/2)}.
\end{split}
\end{equation*}

\hfill$\square$

\section{Proof of Theorem \ref{main_result_theorem}}\label{section_proof_of_the_main_result}

According to \eqref{coefficients_for_orthonormal_generalized_Gegenbauer_polynomials}, \eqref{asymptotic_for_maximum_of_Jacobi_polynomials}, \eqref{values_at_end_points}, and \eqref{asymptotic_for_coefficients_for_orthonormal_generalized_Gegenbauer_polynomials}, we get
\begin{equation*}
\begin{split}
\bigl\|\widetilde{C}_{2n}^{(\lambda,\mu)}\bigr\|_\infty&=\widetilde{a}_{2n}^{\,(\lambda,\mu)} \max\limits_{-1\leq t\leq 1} |P_n^{(\lambda-1/2,\mu-1/2)}(2t^2-1)|=\\
&=\widetilde{a}_{2n}^{\,(\lambda,\mu)} \max\limits_{-1\leq t\leq 1} |P_n^{(\lambda-1/2,\mu-1/2)}(t)|\asymp n^{\max(\lambda,\mu)}.
\end{split}
\end{equation*}
Thus, Theorem \ref{main_result_theorem} is proved for even orthonormal generalized Gegenbauer polynomials.

Consider the case that $\lambda\geq\mu+1$. 

In particular, we have $\lambda>\mu$. By \eqref{asymptotic_for_maximum_of_Jacobi_polynomials}, \eqref{values_at_end_points}, the assumption $\lambda\geq\mu+1$ implies that
\begin{equation*}
\bigl\|P_n^{(\lambda-1/2,\mu+1/2)}\bigr\|_{\infty}=P_n^{(\lambda-1/2,\mu+1/2)}(1)\asymp n^{\lambda-1/2}.
\end{equation*}
Hence, it follows from the above asymptotic equality, \eqref{coefficients_for_orthonormal_generalized_Gegenbauer_polynomials}, and \eqref{asymptotic_for_coefficients_for_orthonormal_generalized_Gegenbauer_polynomials}, that
\begin{equation*}
\bigl\|\widetilde{C}_{2n+1}^{(\lambda,\mu)}\bigr\|_\infty\asymp n^\lambda.
\end{equation*}

Consider the case that $\lambda<\mu+1$.

Note that $2\sin^2\frac{\theta}{2}-1=-\cos\theta$. Making the change of variable $x=\sin\frac{\theta}{2}$ and applying \eqref{coefficients_for_orthonormal_generalized_Gegenbauer_polynomials}\,--\,\eqref{connection_with_plus_minus}, and Lemma \ref{main_result_lemma} with $\alpha=\mu+\frac{1}{2}$, $\beta=\lambda-\frac{1}{2}$, we get
\begin{equation}\label{main_asymptotic_upper_bound}
\begin{split}
\bigl\|\widetilde{C}_{2n+1}^{(\lambda,\mu)}\bigr\|_\infty&=\widetilde{a}_{2n+1}^{\,(\lambda,\mu)}\max\limits_{-1\leq t\leq 1} |t P_n^{(\lambda-1/2,\mu+1/2)}(2t^2-1)|=\\&=\widetilde{a}_{2n+1}^{\,(\lambda,\mu)}\max\limits_{0\leq x\leq 1} |x P_n^{(\lambda-1/2,\mu+1/2)}(2x^2-1)|=\\
&=\widetilde{a}_{2n+1}^{\,(\lambda,\mu)}\max\limits_{0\leq\theta\leq\pi} \Bigl\{\sin\frac{\theta}{2}\,|P_n^{(\lambda-1/2,\mu+1/2)}(-\cos\theta)|\Bigr\}=\\
&=\widetilde{a}_{2n+1}^{\,(\lambda,\mu)}\max\limits_{0\leq\theta\leq\pi} \Bigl\{\sin\frac{\theta}{2}\,|P_n^{(\mu+1/2,\lambda-1/2)}(\cos\theta)|\Bigr\}\lesssim\left\{\begin{array}{lrrll}
n^{\mu},&\mu\!\!\!&>\lambda,&\\[0.4em]
n^{\lambda},&\mu\!\!\!&\leq\lambda.&
\end{array}\right.
\end{split}
\end{equation}

Using \eqref{coefficients_for_orthonormal_generalized_Gegenbauer_polynomials}, \eqref{values_at_end_points}, and \eqref{asymptotic_for_coefficients_for_orthonormal_generalized_Gegenbauer_polynomials}, we obtain
\begin{equation}\label{first_asymptotic_lower_bound}
\bigl\|\widetilde{C}_{2n+1}^{(\lambda,\mu)}\bigr\|_{\infty}\geq \widetilde{C}_{2n+1}^{(\lambda,\mu)}(1)=\widetilde{a}_{2n+1}^{\,(\lambda,\mu)} P_n^{(\lambda-1/2,\mu+1/2)}(1)\gtrsim n^{\lambda}.
\end{equation}
Because of \eqref{main_asymptotic_upper_bound}, \eqref{asymptotic_for_coefficients_for_orthonormal_generalized_Gegenbauer_polynomials}, and \eqref{asymptotic_of_Jacobi_polynomials_on_the_sequence_of_special_points}, we have
\begin{equation}\label{second_asymptotic_lower_bound}
\begin{split}
\bigl\|\widetilde{C}_{2n+1}^{(\lambda,\mu)}\bigr\|_\infty&=\widetilde{a}_{2n+1}^{\,(\lambda,\mu)}\max\limits_{0\leq\theta\leq\pi} \Bigl\{\sin\frac{\theta}{2}\,|P_n^{(\mu+1/2,\lambda-1/2)}(\cos\theta)|\Bigr\}\geq\\
&\geq \widetilde{a}_{2n+1}^{\,(\lambda,\mu)}\sin\frac{n^{-1}}{2}\,|P_n^{(\mu+1/2,\lambda-1/2)}(\cos n^{-1})|\gtrsim n^{\mu}.
\end{split}
\end{equation}

Combining \eqref{main_asymptotic_upper_bound}, \eqref{first_asymptotic_lower_bound}, and \eqref{second_asymptotic_lower_bound}, we observe the desired asymptotic behavior in the considered situation.

Theorem \ref{main_result_theorem} is completely proved.

\section{Conclusion}

The generalized Gegenbauer polynomials play an important role in Dunkl harmonic analysis (see, for example, \cite{dunkl_xu_book_orthogonal_polynomials_2014,dai_xu_book_approximation_theory_2013}). So, the study of these polynomials for different purposes is very natural.

As an application of Theorem \ref{main_result_theorem}, we are going to establish the Hausdorff\,--\,Young and the Hardy\,--\,Littlewood-type inequalities for orthonormal generalized Gegenbauer polynomials in future publications.


\begin{Biblioen}

\bibitem{andrews_askey_roy_book_special_functions_1999}G.\,E.~Andrews, R.~Askey, and R.~Roy, \textit{Special Functions}, Encyclopedia of Mathematics and its Applications \textbf{71}, Cambridge University Press, Cambridge, 1999.

\bibitem{dai_xu_book_approximation_theory_2013}F.~Dai and Y.~Xu, \textit{Approximation theory and harmonic analysis on spheres and balls}, Springer Monographs in Mathematics, Springer, 2013.

\bibitem{dunkl_xu_book_orthogonal_polynomials_2014}C.\,F.~Dunkl and Y.~Xu, \textit{Orthogonal polynomials of several variables}, 2nd ed., Encyclopedia of Mathematics and its Applications \textbf{155}, Cambridge University Press, Cambridge, 2014.

\bibitem{osilenker_book_orthogonal_polynomials_1999}B.~Osilenker, \textit{Fourier series in orthogonal polynomials}, World Scientific, Singapore, 1999.

\bibitem{szego_book_orthogonal_polynomials_1975}G.~Szeg\"{o}, \textit{Orthogonal polynomials}, 4th ed., American Mathematical Society Colloquium Publications \textbf{23}, American Mathematical Society, Providence, Rhode Island, 1975.

\end{Biblioen}

\noindent \textsc{Department of scientific research, Tula State University, Tula, Russia }

\noindent \textit{E-mail address}: \textbf{veprintsevroma@gmail.com}

\end{document}